\documentclass{amsart}

\usepackage{amsmath,amssymb,amsthm,dsfont,hyperref,mathtools,wasysym}
\usepackage[active]{srcltx}
\usepackage[top=3.2cm,bottom=3.2cm,left=3.2cm,right=3.2cm]{geometry}

\usepackage{fancyhdr}
\usepackage{enumitem}
\usepackage{trfsigns}
\usepackage{graphicx}
\usepackage{xcolor}

\usepackage{upref} 

\newcommand{\divo}{\textnormal{div}}
\newcommand{\Log}{\textnormal{L$^2$(Log L)$^\alpha$}}

\def\en{\mathbb N}

\def\er{\mathbb R}
\def\eqn#1$$#2$${\begin{equation}\label#1#2\end{equation}}

\newcommand\dx{\,dx}

\newcommand\eps\varepsilon
\renewcommand\rho\varrho
\newcommand\g\gamma
\newcommand\G\Gamma

\renewcommand\L\Lambda
\newcommand\loc{{\rm loc}}

\newcommand\qq\qquad
\newcommand\R{\mathds{R}}

\newcommand\vp\varphi
\newcommand\ts\textstyle

\newcommand{\mint}{\mathop{\int\hskip -1,05em -\, \!\!\!}\nolimits}

\newcommand{\kabs}[1]{\ensuremath{\vert#1\vert}}

\newtagform{boldbrackets}{\hspace{8ex}\bf(}{\bf)}
\newtagform{standardbrackets}{)}{)}



\newcommand{\comment}[1]{}

\allowdisplaybreaks[3]

\newtheorem{theorem}{Theorem}[section]

\usepackage[normalem]{ulem}
\newcommand{\stkout}[1]{\ifmmode\text{\sout{\ensuremath{#1}}}\else\sout{#1}\fi}

\begin{document}

\title[Lipschitz criteria for non-uniformly elliptic problems]{Optimal Lipschitz criteria and local estimates for non-uniformly elliptic problems}

\author[L. Beck]{Lisa Beck}
\address{Institut f\"{u}r Mathematik, Universit\"{a}t Augsburg, Universit\"{a}tsstr. 14, 86159~Augsburg, Germany}
\email{lisa.beck@math.uni-augsburg.de}

\author[G. Mingione]{Giuseppe Mingione}
\address{Dipartimento SMFI, Universit\`a di Parma, Viale delle Scienze 53/a, Campus, 43124 Parma, Italy}
\email{giuseppe.mingione@unipr.it}



\maketitle
\vspace{3mm}

\centerline{To Vladimir Gilelevich Maz'ya,} \centerline{master of analysis across two different worlds and centuries,}
\centerline{on his 80th birthday}
\begin{abstract} We report on new techniques and results in the regularity theory of general non-uniformly elliptic variational integrals. By means of a new potential theoretic approach we reproduce, in the non-uniformly elliptic setting, the optimal criteria for Lipschitz continuity known in the uniformly elliptic one and provide a unified approach between non-uniformly and uniformly elliptic problems. \end{abstract}

\section{Introduction and results}
\label{sec-intro}
Non-uniformly elliptic equations and functionals are a classical topic in partial differential equations and in the Calculus of Variations. They emerge in many different contexts, and they are often stemming from geometric and physical problems \cite{DE, LU, Simon, UU, Zhikov}. The study of regularity of their solutions involves a wealth of methods and techniques. Here we announce a few results from~\cite{BM}, concerning the regularity of minimizers of non-uniformly elliptic variational integrals of the type
\eqn{genF}
$$
W^{1,1}(\Omega;\er^N) \ni w \mapsto \mathcal F(w;\Omega)\coloneqq  \int_{\Omega}  \left[ F(Dw)-f\cdot w \right]\, dx\;,   
$$ 
where $ \Omega \subset \er^n $ is an open subset, $n\geq 2$, $N\geq 1$, and 
the integrand $F\colon \er^{N\times n}\to \er$ satisfies suitable convexity and growth conditions.  

The aim of~\cite{BM} is twofold, specifically:
\begin{itemize}
\item To identify sharp conditions on the datum $f$ in order to guarantee that minimizers are locally Lipschitz continuous.  
\item To introduce a new potential theoretic technique allowing to reduce the treatment of non-uniformly elliptic problems to the one of uniformly elliptic ones. This technique yields optimal and new local estimates already in the case $f\equiv 0$. 
\end{itemize}
Let us immediately clarify the meaning or the terminology concerning the non-uniform ellipticity. 
The Euler--Lagrange equation of the functional in~\eqref{genF} formally reads as 
$$
-\divo\, \partial F(Du)=f\;. 
$$
Ellipticity is then tested by looking at the so-called second variation equation, and therefore 
at the ellipticity ratio
\eqn{ratio}
$$
\mathcal R(z) \coloneqq \frac{\mbox{highest eigenvalue of}\ \partial^2  F(z)}{\mbox{lowest eigenvalue of}\  \partial^2 F(z)} \;.
$$
The term non-uniform ellipticity refers to the situation when $\mathcal R(z) $ is not bounded. This is the case of interest here and therefore is not present in more classical situations as for instance the one of the $p$-Laplacean operator $F(z)=:|z|^p/p$, for $p>1$ (see \cite{Manth, Uh, Ur} for basic regularity results). More in general, when considering functionals of the type 
\eqn{doppio}
$$
w \mapsto \int_{\Omega}  \left[A(|Dw|)-f\cdot w \right]\, dx\;, \qquad \text{where } A(t)\coloneqq  \int_0^t \tilde a(s)s\, ds\;,
$$
then the non-uniform ellipticity is connected to the absence of the double-side bound
\eqn{uni}
$$
-1< i_a \leq \frac{\tilde a'(t)t}{\tilde a(t)}\leq s_a< \infty
$$
(considered e.g.~in~\cite{Baroni, CM0, CM1, L1}). This implies the validity of the so-called $\Delta_2$-condition
\eqn{delta}
$$
A(2t)\lesssim A(t)\;, 
$$
and therefore immediately allows for polynomial controls for $t\mapsto A(t)$. 
In the classical uniformly elliptic case \eqref{doppio}-\eqref{uni}, the question of determining sharp conditions on $f$ implying the local and global Lipschitz continuity of minimizers has found a satisfying solution only in the last few years. Results are available both for the scalar and for the vectorial case, see e.g.~\cite{CM1, CM3, DM, KUUMIN12j, KUUMIN13, KUUMINpress, mis3}. The most general answer can be obtained in terms of Riesz potentials, and we refer to~\cite{KUUMIN14} for an overview of nonlinear potential estimates. In terms of function spaces, a neat and sharp condition on $f$ makes use of Lorentz spaces  
$L(n,1)(\Omega)$, that is
\eqn{cond-lo}
$$
f \in L(n,1)(\Omega) :\Longleftrightarrow \| f\|_{L(n,1)(\Omega)}\coloneqq  \int_0^\infty |\{x \in \Omega \colon |f(x)|> \lambda\}|^{1/n} \, d\lambda <\infty\;.
$$
The remarkable feature is that condition~\eqref{cond-lo} does not depend on the functional under consideration. In particular, when considering the $p$-Laplacean equation $
-\divo (|Du|^{p-2}Du)=f$, i.e., for $F(z)=|z|^p/p$ in~\eqref{genF}, the condition to ensure the local Lipschitz continuity of minimizers is in fact independent of $p$, as first observed in~\cite{DM}. In particular, condition~\eqref{cond-lo} is sharp for the Poisson equation $-\triangle u=f$, as a consequence of a well-known theorem of Stein~\cite{St} and the work of Cianchi \cite{CiaG}. Notice that the Lorentz space $L(n,1)$ is an intermediate, borderline space in the sense that $L^{n+\eps}\subset L(n,1)\subset L^n$ holds for every $\eps>0$ (and all such inclusions are actually strict). Concerning the non-uniformly elliptic case, essentially nothing is known. 

A main outcome of~\cite{BM} is that the effectiveness of condition~\eqref{cond-lo} is extended to the case 
of non-uniformly elliptic problems, thereby filling a huge gap in the literature. We also notice that, when looking at non-uniformly elliptic functionals, the results in~\cite{BM} are new already in the case when $f$ is more regular than in~\eqref{cond-lo}. The usual cases treated in the literature are indeed $f\equiv 0$ or $f\in L^\infty$, where the latter one was only considered under polynomial growth conditions (see~\cite{M91}). As we shall see in Section~\ref{unisec} below, the techniques of~\cite{BM} yield new results in the classical uniformly elliptic case with \eqref{doppio}-\eqref{uni} as well. The class of functionals treated in~\cite{BM} is large, and essentially covers all the main models appearing in the literature as far as autonomous functionals are considered. It includes functionals with polynomial, yet unbalanced growth conditions, i.e., when $F$ satisfies different growth and coercivity conditions such as 
\eqn{pqpq} 
$$
|z|^p \lesssim F(z) \lesssim |z|^q+1\,, \qquad 1< p\leq q\;.
$$
These are usually called $(p,q)$-growth conditions, according to Marcellini's terminology from~\cite{M91}. We also cover the case of functionals with very fast growth as for instance 
\eqn{superexp}
$$
w \mapsto \int_{\Omega} \Big[\exp(\exp(\dots \exp(|Dw|^p)\ldots))-f\cdot w \Big]\, dx \;,  \quad p \geq 1\;,$$
and for which the $\Delta_2$-condition~\eqref{delta} fails.

In the following, we report a few results from~\cite{BM}, essentially referring to the two main model situations~\eqref{pqpq} and~\eqref{superexp}. These are in turn consequences of a more general result still from~\cite{BM} whose formulation turns out to be rather technical as it is devised to cover large and different classes of model cases. For a reasonably updated review of the existing regularity literature on the subject we refer to~\cite{Dark}.

For the rest of the paper, in view of condition~\eqref{cond-lo}, we assume that $f\in L^n(\Omega;\R^N)$. Moreover, we recall that a function $u \in W^{1,1}_{\loc}(\Omega;\er^N)$ is a \emph{local minimizer} of the functional $\mathcal F$ in~\eqref{genF} with $f \in L^n_{\loc}(\Omega;\R^N)$ if, for every open subset $\tilde \Omega\Subset \Omega$, we have  that $\mathcal F(u;\tilde \Omega) <\infty$ and moreover $\mathcal F(u;\tilde \Omega)\leq \mathcal F(w;\tilde \Omega)$ holds whenever $w \in u + W^{1,1}_0(\tilde \Omega; \er^N)$. As $f \in L^n_{\loc}(\Omega;\R^N)$, Sobolev embedding theorem implies $f\cdot u \in L^1_{\loc}(\Omega)$ so that $F(Du)\in L^1_{\loc}(\Omega)$. In the following we denote by $ B_r(x_0)\coloneqq \{x \in \er^n \colon |x-x_0|< r\}$ the open ball with center $x_0$ and radius $r>0$; the center will be omitted ($B_r \equiv B_r(x_0)$) when clear from the context. If not otherwise stated, different balls in the same context will share the same center. With $\mathcal B \subset \er^{n}$ being a measurable subset with bounded positive measure $0<|\mathcal B|<\infty$, and with $g \colon \mathcal B \to \er^{k}$, $k\geq 1$, being a measurable map, we shall denote its integral average by  
$$
   (g)_{\mathcal B} \equiv \mint_{\mathcal B}  g(x) \, dx  \coloneqq  \frac{1}{|\mathcal B|}\int_{\mathcal B}  g(x) \, dx\;.
$$

\section{A model result for polynomial growth conditions}\label{pqsec} 
In this section we give a model result covering the case of functionals under $(p,q)$-growth conditions as in~\eqref{pqpq}. They have been studied in detail in the literature and, in the setting of the Calculus of Variations, an extensive treatment has been provided by Marcellini~\cite{M91}; see also~\cite{schmidt} for the vectorial case. In this setting, the integrand $F \colon \er^n \to \er$ is convex function and locally $C^2$-regular in $\er^n\setminus\{0\}$. The crucial growth and ellipticity properties of $F$ are described as follows:
\eqn{asp1}
$$
\left\{
\begin{array}{c}
\nu (|z|^2+\mu^2)^{p/2} \leq F(z) \leq L(|z|^2+\mu^2)^{q/2}+L(|z|^2+\mu^2)^{p/2}\\[5 pt]
(|z|^2+\mu^2) |\partial^2 F(z)| \leq L(|z|^2+\mu^2)^{q/2}+L(|z|^2+\mu^2)^{p/2}\\ [5 pt]
\nu (|z|^2+\mu^2)^{(p-2)/2}|\xi|^2\leq  \langle \partial^ 2F(z)\xi, \xi\rangle\;, 
 \end{array}\right.
$$
for every choice of $z, \xi \in \er^n$ with $z\not=0$ and for exponents $1 \leq p \leq q$. Here $0<\nu\leq 1 \leq L$ are fixed ellipticity constants and $\mu \in [0,1]$ serves  to distinguish the so-called degenerate case ($\mu=0$) and non-degenerate case ($\mu >0$). In this case the ellipticity ratio in~\eqref{ratio} is controlled as follows:
\eqn{rateratio}
$$
\mathcal R(z) \,
\lesssim  (z|^2+\mu^2)^{(q-p)/2} +1\;. 
$$
It therefore becomes unbounded as $|z|\to \infty$ for the gradient variable, with a speed which is proportional to the so-called gap $q-p$. The main result for $(p,q)$-growth functionals is the following

\begin{theorem} [Scalar $(p,q)$-estimates]\label{main1} Let $u \in W^{1,1}_{\loc}(\Omega)$ be a local minimizer of the functional $\mathcal F$ in~\eqref{genF} under assumptions~\eqref{asp1} with $1 < p \leq q$ and $n>2$. Assume 
\eqn{ass1}
$$
\frac qp < 1+ \min\left\{ \frac 2n, \frac{4(p-1)}{p(n-2)}\right\}\qquad \mbox{and}\qquad f \in L(n,1)(\Omega) \;.
$$
Then $Du$ is locally bounded in $\Omega$. Moreover, the local a 
priori estimate
\begin{multline}
\|Du\|_{L^{\infty}\left(B/2\right)} \leq c\left(\mint_{B}F(Du)\dx +\| f\|_{L(n,1)(B)}^{\frac{p}{p-1}}\right)^{\frac 1p}
\\
 \ \ + c \left(\mint_{B}F(Du)\dx + \| f\|_{L(n,1)(B)}^{\frac{p}{p-1}} \right)^{\frac{2}{(n+2)p-nq}}
+c\left[\| f\|_{L(n,1)(B)}\right]^{\frac{4}{4(p-1)-(n-2)(q-p)}}  
\label{estimatepq}
\end{multline}
holds for a constant $c\equiv c(n,p,q,\nu, L)$, whenever $B\Subset \Omega$ is a ball. When $p\geq 2-4/(n+2)$ or when $f\equiv 0$, condition~\eqref{ass1} can be replaced by 
\eqn{marbound} 
$$
\frac qp < 1+\frac 2n\;.
$$
\end{theorem}

The above result is valid for scalar minimizers. When passing to the vectorial case, singularities may in general emerge, cf.~\cite{Dark}. Minimizers are regular outside so-called singular sets whose dimension can be estimated see e.g.~\cite{KrM2, schmidt}; so called partial regularity theory comes into the play. However, assuming a suitable, quasidiagonal structure, one is able to recover everywhere regularity together with the scalar results. This is known since the fundamental work of Uhlenbeck~\cite{Uh}; as for the non-homogeneous case, we refer to the recent paper~\cite{KUUMINpress} that also contains potential estimates. In our situation we indeed have 

\begin{theorem}[Vectorial $(p,q)$-estimates]\label{main1-vec} Let $u \in W^{1,1}_{\loc}(\Omega;\er^N)$ be a local minimizer of the functional $\mathcal F$ in~\eqref{genF} 
under assumptions~\eqref{asp1} with $1 < p \leq q$ and $n>2$. Assume 
\eqn{ass1-vec}
$$
\frac qp < 1+ \min\left\{ \frac 1n, \frac{2(p-1)}{p(n-2)}\right\}\qquad \mbox{and}\qquad f \in L(n,1)(\Omega) \;,
$$
and that there exists a 
$C^1_{\loc}[0,\infty)\cap C^2_{\loc}(0,\infty)$-regular function $\tilde F \colon [0, \infty) \to [0, \infty)$ such that $F(z)=\tilde F(|z|)$ for every $z \in \er^{N\times n}.$ Finally, assume that the function $t \mapsto (\mu^2+t^2)^{(2-\gamma)/2}\tilde F'(t)/t$ is non-decreasing for some $\gamma>1$; in particular, $\gamma=p$ is admissible. 
Then $Du$ is locally bounded in $\Omega$ and an estimate similar to~\eqref{estimatepq} holds. In the case that $t \mapsto \tilde F'(t)/t$ is non-decreasing, assumption~\eqref{ass1-vec} can be relaxed to~\eqref{ass1}. 
\end{theorem}

Conditions~\eqref{ass1},~\eqref{marbound} and~\eqref{ass1-vec} are standard in the present setting and serve to bound the rate of possible blow-up of $\mathcal R(z)$ as described in~\eqref{rateratio}. They occur in various forms, cp.~\cite{BCM, BCM2, BiFu, CKP, COLMIN15o, COLMIN15, DF, L3, M91}, and counterexamples show that they are in fact necessary~\cite{M91}. In particular,~\eqref{marbound} is assumed in several papers starting from~\cite{M91}.  

Estimate~\eqref{ass1} is general enough to reproduce several classical results when applied to particular settings. Indeed, when dealing with the $p$-Laplacean equation $
-\divo (|Du|^{p-2}Du)=f$, we then take $p=q$ and $\mu=0$, and~\eqref{estimatepq} gives the local estimate (see for instance \cite{DM, KUUMIN13, KUUMINpress})
\eqn{optimal}
 $$
 \|Du\|_{L^{\infty}\left(B/2\right)}\lesssim  \left(\mint_{B} |Du|^p \dx\right)^{\frac 1p}+\|f\|_{L(n,1)(B)}^{\frac{1}{p-1}}\;. 
 $$
This is nothing but the classical $L^\infty-L^p$ estimate for $p$-harmonic functions~\cite{Manth} when $f\equiv 0$. Instead, keeping $p\not=q$ as in~\eqref{ass1}, but considering the homogeneous case $f\equiv 0$, it gives back the following classical bound of Marcellini \cite[Theorem 3.1]{M91} (there $\mu=1$, $q \geq p\geq 2$ and~\eqref{marbound} are assumed):
$$
\|Du\|_{L^{\infty}\left(B/2\right)} \lesssim   \left(\mint_{B}F(Du)\dx\right)^{\frac{2}{(n+2)p-nq}} +1 \;. 
$$

Finally, let us comment on the two dimensional case $n=2$. For this we use a different, almost sharp characterization making use of suitable Orlicz spaces. We recall that, given a Young function $A\colon [0, \infty)\to[0, \infty)$, with $\Omega\subset \er^n$, the Orlicz space $L^A(\Omega)$ is defined as the vector space of measurable maps $g$ such that $\|A(|g|)\|_{L^1(\Omega)} < \infty$ the following Luxemburg norm is finite:
$$
\|w\|_{L^A(\Omega)}\coloneqq  \inf\left\{ \lambda >0 \, \colon \, \int_{\Omega} A\left(\frac{|w|}{\lambda}\right)\, dx \leq 1\right\}\;.
$$
This is a Banach space. The case of interest for us is given by the choice $A(t)= t^2\log^\alpha(\e + t)$ for $\alpha \geq 1$, for which the standard notation is $L^A\equiv \Log$. For the Poisson equation $-\triangle u\in \Log$ it is known that $\alpha>1$ is sufficient for local Lipschitz continuity of $u$, and this result is sharp in this scale of spaces. In~\cite{BM} we are able to reproduce almost the same criterion; indeed it holds the following:
\begin{theorem}[The two dimensional case]\label{main1-two} Let $u \in W^{1,1}_{\loc}(\Omega)$ be a local minimizer of the functional in~\eqref{genF} under assumptions~\eqref{asp1} with $1 < p \leq q$ and $n=2$. Assume
\begin{equation*}
q<2p \qquad \mbox{and}\qquad f \in \Log(\Omega) \mbox{ for some $\alpha>2$,} 
\end{equation*}
Then $Du$ is locally bounded in $\Omega$. \end{theorem} 

The previous result comes along with local estimates and extends to the vectorial case (under suitable structure conditions mentioned above) as well; see~\cite{BM}. 

\section{Non-polynomial growth conditions and natural estimates}\label{expsec} 
Here we deal with non-polynomial growth conditions, and we present the results directly in the vectorial case. An important example of a functional with fast growth conditions appears in~\eqref{superexp}. It belongs to a more general family which is described as follows. 
Let $\{p_k\}$ be a sequence of 
real numbers such that $p_0 > 1$, $p_k > 0$ for every $k\in \en$. By inductions, we define ${\bf e}_{p_k}\colon [ 0, \infty) \to \er$ as
$$\left\{
 \begin{array}{ccc}
 {\bf e}_{k+1}(t) &\coloneqq  & \exp \left[\left({\bf e}_{k}(t)\right)^{p_{k+1}}\right] \\ [4 pt]{\bf e}_{0}(t) &\coloneqq & \exp(t^{p_0})\;.
 \end{array}
 \right.
$$
and consider the variational integrals
\eqn{espik}
$$
W^{1,1}(\Omega;\er^N) \ni w \mapsto \mathcal E_k(w)\coloneqq   \int_{\Omega}  \left[ {\bf e}_{k}(|Dw|)-f w \right]\, dx\;.
$$
We then have
\begin{theorem}[Exponential estimates]\label{main4} Let $k\geq 0$ be an integer and let $u \in W^{1,1}_{\loc}(\Omega; \er^N)$ be a local minimizer of the functional $\mathcal E_k$ in~\eqref{espik}. 
\begin{itemize}
\item If $f \in L(n,1)(\Omega;\er^{N})$ and $n>2$, then $Du$ is locally bounded in $\Omega$.
\item When $n=2$ the same conclusion holds provided $f \in \Log(\Omega)$ for some $\alpha>2$.
\item  Finally, when $f\equiv 0$, the local estimate 
\eqn{thecasegeneral}
$$
\|Du\|_{L^{\infty}\left(B/2\right)} \leq c \, {\bf e}_{k}^{-1}\left(\mint_{B}  {\bf e}_{k}(|Du|) \dx\right)+c
$$
holds for a constant $c\equiv c(n,N,k,p_0, \ldots, p_k)$ and for every ball $B \Subset \Omega$. 
\end{itemize}
\end{theorem}
The above theorem features new and relevant information already in the case $f\equiv 0$ (originally treated by Marcellini~\cite{M96, M96b}; see also~\cite{DE, L2}). To understand the progress, let us consider the simplest example given by $$w \mapsto \int_{\Omega} \exp(|Dw|^p)\, dx\;,\quad p>1\;.$$
For this case the best a priori estimate available up to now, obtained in~\cite{M96, M96b}, reads as
\eqn{thecase}
$$ \|Du\|_{L^{\infty}\left(B/2\right)} \leq c_{\eps}\left(\mint_{B} \exp (|Du|^p) \dx \right)^{1+\eps} +c_{\eps} \;, \quad \mbox{for every}\ \eps >0\;,$$
whereas our estimate~\eqref{thecasegeneral}, when applied in this situation, gives
$$
 \|Du\|_{L^{\infty}\left(B/2\right)}^p\leq c \log \left(\mint_{B} \exp (|Du|^p) \dx\right)+c\;,\qquad c\equiv c(n,N,p)\;. 
$$
This estimate is new and parallels the one for $p$-harmonic functions in~\eqref{optimal} in the sense that it exhibits a gain of an exponential scale with respect to~\eqref{thecase}. The gain with respect to the estimates in~\cite{M96, M96b} increases when considering functionals of even faster growth as in~\eqref{superexp}. Estimate~\eqref{thecasegeneral} is a special occurrence of a more general result obtained in~\cite{BM} which concerns the special situation of functionals defined in Orlicz spaces of the type
\eqn{orlicz}
$$
w \mapsto  \int_{\Omega} A(|Dw|)\, dx\;,
$$
where $A\colon [0, \infty) \to [0, \infty)$ is a Young function. The case $A(t)=t^p/p$ for $p\geq 1$ 
is the most basic example. A natural question is of course whether or not natural a priori estimates of the type
$$
\|Du\|_{X}\lesssim A^{-1} \left(\mint_{B} A(|Du|) \dx\right)+ 1
$$
hold for function spaces $X$ that embed into $L^A(B/2)$. The answer is known to be positive \cite{Baroni, CF, FS} provided the $\triangle_2$-condition~\eqref{delta} holds. The techniques in~\cite{BM} allow to give a general positive answer to this question also in the more difficult situation when the $\triangle_2$-condition is dropped. The results in~\cite{BM} cover the most important case of $X\equiv L^\infty(B/2)$ and a large class of functionals of the type~\eqref{orlicz}. 

\section{Revisiting uniformly elliptic operators}\label{unisec}
The results in~\cite{BM} allow to get new conclusions in the classical uniformly elliptic case too. We indeed have the following:
\begin{theorem}[Natural growth estimates under $\triangle_2$-condition]\label{maing-uniform} 
Let $u \in W^{1,1}_{\loc}(\Omega; \er^N)$ be a local minimizer of the functional in~\eqref{doppio}, with
$A \colon [0, \infty) \to [0, \infty)$ being of class $C^1_{\loc}[0,\infty)\cap C^2_{\loc}(0,\infty)$. Assume that the uniformly ellipticity assumptions~\eqref{uni} hold and that $f \in L(n,1)(\Omega;\er^N)$ for $n>2$. Then $Du\in L^{\infty}_{\loc}(\Omega;\er^{N\times n})$ and the estimate
\eqn{esti-unifo}
$$
\|Du\|_{L^{\infty}(B/2)}\leq cA^{-1}\left( \mint_{B} A(|Du|)\dx \right) + c
A^{-1}\left(\| f\|_{L(n,1)(B)}^{\frac{i_{a}+2}{i_{a}+1}} \right) + c_2
$$
holds for every ball $B \Subset \Omega$, for constants $c, c_2$ depending only on $n,N,i_a,s_a$ and $\tilde a (1)$. When $n=2$ a similar results holds assuming that $f \in \Log(\Omega;\er^N)$ for some $\alpha>2$. In the case it is
$$i_2\coloneqq  \liminf_{t\to 0}\,\frac{\tilde a (t)}{t^{i_a}}>0\;,$$ then we can take $c_2=0$ in~\eqref{esti-unifo} and the constant $c$ depends also on $i_2$. 
\end{theorem}
The above result extends those of Baroni~\cite{Baroni} and Lieberman~\cite{L1} when the problem is vectorial. It also provides a local analog to the global bounds of Cianchi \& Maz'ya \cite{CM0, CM1, CM3}. It is worth remarking that, when applied to the usual $p$-Laplacean case $A(|z|)\equiv |z|^p/p$, estimate~\eqref{esti-unifo} reduces to the classical one in~\eqref{optimal}.  

\section{Existence and regularity for equations under $(p,q)$-growth conditions}
When considering general equations under $(p,q)$-growth conditions, not necessarily arising from variational integrals, it is still possible to prove existence of locally Lipschitz continuous solutions. Specifically, let us consider Dirichlet problems of the type
\eqn{dir1}
$$
\left\{
\begin{array}{c}
-\divo\, a (Du)=f \quad \mbox{in}\ \Omega\\ [3 pt]
u \equiv u_0  \quad \mbox{on}\ \partial \Omega\;,
 \end{array}\right.
\qquad \qquad 
u_0 \in W^{1,\frac{p(q-1)}{p-1}}(\Omega)\;.
$$
Here we additionally assume that $\Omega$ is a bounded and Lipschitz regular domain. The vector field $a\colon \er^n \mapsto \er^n$ is assumed to be $C^1$-regular outside the origin and such that
\eqn{asp1-ex}
$$
\left\{
\begin{array}{c}
|a(z)|+(|z|^2+\mu^2)^{\frac 12} |\partial a(z)| \leq L(|z|^2+\mu^2)^{\frac{q-1}{2}}+L(|z|^2+\mu^2)^{\frac{p-1}{2}}\\ [5 pt]
\nu (|z|^2+\mu^2)^{\frac{p-2}{2}}|\xi|^2\leq  \langle \partial a(z)\xi, \xi\rangle\;
 \end{array}\right.
$$
holds for every choice of $z, \xi \in \er^n$ with $z\not=0$, and for exponents $1 \leq p \leq q$, ellipticity constants $0<\nu\leq 1 \leq L$ and $0 \leq \mu \leq 1$. Marcellini's by now classical result~\cite{M91} states the existence of locally Lipschitz continuous solutions to~\eqref{asp1-ex} assuming that $f \in L^\infty(\Omega)$. Marcellini's result is upgraded in~\cite{BM} to the sharp level $f \in L(n,1)$. The precise statement of the result is as follows:

\begin{theorem}[Existence of locally Lipschitz solutions]\label{main-eq} With 
$1 <  p \leq q$, $n>2$,  and under assumptions~\eqref{ass1-vec} and~\eqref{asp1-ex}, 
there exists a solution $u \in W^{1, \infty}_{\loc}(\Omega)\cap W^{1,p}(\Omega)$ to the Dirichlet problem in~\eqref{dir1}. Moreover,  the a priori estimate
\eqn{estimatepq-eq}
$$ \|Du\|_{L^{\infty}\left(B/2\right)}  \leq   c\frac{\mathcal D}{|B|}+c\left(\frac{\mathcal D}{|B|}\right)^{\frac{p}{p-(q-p)n}}
+ c\|f\|_{L(n,1)(B)}^{\frac{1}{p-1}}+c\|f\|_{L(n,1)(B)}^{\frac{2}{2(p-1)-(q-p)(n-2)}}
$$
holds for every ball $B\Subset \Omega$, where $c\equiv c (n,p,q,\nu, \Lambda)$, and 
$$ \mathcal D^p\coloneqq  \int_{\Omega}\left(|Du_0|^2+1\right)^{\frac{p(q-1)}{2(p-1)}}\, dx +  c|\Omega|\|f\|_{L^n(\Omega)}^{p/(p-1)}\;.$$
 When $n=2$ and $f \in \Log(\Omega)$ holds for some $\alpha>2$, again $Du$ is locally bounded in $\Omega$ and an estimate similar to~\eqref{estimatepq-eq} holds upon replacing
$\| f\|_{L(n,1)(B)}$ by $\|f\|_{\Log(B)}$. 
\end{theorem}
We notice that, in the case when $p=q$ and $\Omega$ coincides with a ball $B$, estimate~\eqref{estimatepq-eq} gives
$$
\|Du\|_{L^\infty(B/2)}^p\lesssim \mint_{B}\left(|Du_0|^p+1\right)\, dx +  \|f\|_{L(n,1)(B)}^{p/(p-1)}\;,$$
which is the usual a priori local estimate for the Dirichlet problems; see~\cite{KUUMINpress} and compare with~\eqref{optimal}. Apart from the optimal regularity $f \in L(n,1)(\Omega)$ of the right-hand side, further differences compared with~\cite{M91} rely in the fact that we are also able to treat the degenerate case $\mu=0$, while we are not assuming any control on the antisymmetric part of $\partial a (\cdot)$. Specifically, we are assuming no inequality of the type $|\partial a(z) - [\partial a(z)]^{t}| \lesssim |z|^{(q+p-2)/2},$ as indeed done in~\cite{M91}. This is essentially reflected in the more restrictive bound~\eqref{ass1-vec} on $q/p$ with respect to the one considered in~\cite{M91}, which is in fact ~\eqref{marbound}.  

\section{Potential theoretic techniques} At the core of the new approach developed in~\cite{BM} there lies a set of nonlinear potential theoretic techniques that find their origins in the seminal papers of Maz'ya and Havin~\cite{maz, MH}. Potential estimates for solutions to nonlinear partial differential equations have first been developed by  Kilpel\"ainen \& Mal\'y~\cite{KILMAL94}, in turn relying on the fundamental methods of De Giorgi~\cite{DG}. Gradient potential estimates have been derived in \cite{Baroni, KUUMIN13, KUUMIN14} in the scalar case and in~\cite{KUUMINpress} for the vectorial one. We refer to~\cite{KUUMIN12j} for an overview on potential estimates. Such techniques have been successfully applied to the uniformly elliptic case, whilst a general extension to the non-uniformly elliptic one has always appeared problematic. The main novelty of~\cite{BM} is indeed to find a way to use a potential theoretic approach in the non-uniformly elliptic setting. This involves the use of a modified Riesz potential of the right-hand side data $f$, originally introduced in~\cite{DM}, defined by
$$
 {\bf P}_1^f(x_0,R) \coloneqq \int_0^R \left( \rho^{2} \mint_{B_\rho(x_0)} \kabs{f}^2 \dx \right)^{1/2} \frac{d \rho}{\rho} \;.
$$
for $x_0 \in \Omega$ and $R>0$. The potential $ {\bf P}_1^f(\cdot,R)$ plays in the present context the role of the (truncated) Riesz potential which is instead defined by
$$
{\bf I}_{1}^f(x_0,R)\coloneqq  \int_0^R \frac{|f|(B_{\varrho}(x_0))}{\varrho^{n-1}}\, \frac{d\varrho}{\varrho}= 
 \int_0^R \frac{1}{\varrho^n}\int_{B_{\varrho}(x_0)}|f|\, dx \, d\varrho
\;,
$$
provided the argument function $f$ is at least locally $L^1$-regular. This can be easily seen by using H\"older inequality to estimate $
{\bf I}_{1}^f(x_0,R)  \leq|B_1| {\bf P}_1^f(x_0,R)$ and observing that the two quantities share the same homogeneity and scaling properties. In~\cite{BM} we prove an estimate that allows to bound, locally, the $L^\infty$-norm of $Du$ 
in terms of the $L^\infty$-norm of ${\bf P}_1^f(\cdot,R)$ (modulo controllable terms). The catch with Lorentz spaces and local estimates then comes from the relation $$ \|{\bf P}_1^{f}(\cdot,R)\|_{L^{\infty}(B)}\lesssim \| f\|_{L(n,1)(B_{2R})}$$ that holds for every ball $B\subset \er^n$. Similar relations hold in the two dimensional case for the space $\Log$ with $\alpha >2$. The a priori estimates are then combined with an approximation argument aimed at by-passing the fact that in~\cite{BM} we never assume the validity of the $\triangle_2$-condition~\eqref{delta}. Actually different approximations of the functional are needed in the scalar and in the vectorial case. The fact that the we are not assuming the $\triangle_2$-condition~\eqref{delta} poses additional technical difficulties also when building the approximation argument.  
\section{More on non-uniformly elliptic problems}

The extension of results as those in Theorem~\ref{main1} to the case of general non-autonomous functionals of the type
$$
w \mapsto \mathcal F(w;\Omega)\coloneqq  \int_{\Omega}  \left[H(x,Dw)-f\cdot w \right]\, dx 
$$ 
is not an easy task and it is the object of ongoing investigation~\cite{DeMi}. Examples of energies of this type are obviously given by functionals with coefficients of the type 
$$
w \mapsto \mathcal F(w;\Omega)\coloneqq  \int_{\Omega}  \left[c(x)F(Dw)-f\cdot w \right]\, dx\,, \qquad 0 < \nu \leq c(x)\leq L\;,
$$
where $F(\cdot)$ is of the type considered for instance in Theorem~\ref{main1}, or it is one of the exponential growth. Intermediate classes of non-autonomous functionals can be considered as well. These are somehow more interesting as they present a special form of non-uniform ellipticity that can be detected only in a non-local fashion. Indeed, while considering the ratio in~\eqref{ratio} gives a bounded quantity, the non-local ratio
$$
\mathcal R(z, B) \coloneqq \frac{\sup_{x \in B}\,  \mbox{highest eigenvalue of}\ \partial^2  H(x,z)}{\inf_{x \in B}\,\mbox{lowest eigenvalue of}\  \partial^2 H(x,z)} \;, 
$$
does not. This is for instance the case of energies of the type
\eqn{HHaa} 
$$
\left\{
\begin{array}{c}
\displaystyle w \mapsto \mathcal F(w;\Omega)\coloneqq  \int_{\Omega}  \left[|Dw|^{p}+a(x) |Dw|^{q}-f\cdot w \right]\, dx\\ 
[14 pt]
1<p\leq q\,, \quad 0 \leq a(x)\leq L\;,
\end{array}
\right.
$$
as the coefficient $a(x)$ is allowed to be zero. For such energies, ad hoc methods can be developed; see for instance \cite{BCM, BCY, COLMIN15, COLMIN15o, COLMIN16, DeOh}. Another class of relevant non-autonomous functionals is the one involving a variable growth exponent. In this case the relevant energy is given by
$$
w \mapsto \int_{\Omega}  \left[|Dw|^{p(x)}-f\cdot w \right]\, dx\,
$$
and requires less stringent assumptions on the exponent $p(x)$ than those needed on the coefficient $a(x)$ in~\eqref{HHaa}. The regularity problems have been in this case the object of intensive investigation. See the survey~\cite{Dark} for a rapid review and also~\cite{RR} for more general results. Recently, manifold constrained problems have been investigated in the case of variable growth exponent functionals~\cite{DF2} too. In some cases, under more restrictive assumptions, also constrained $(p,q)$-growth conditions can be considered~\cite{DF}.

\bibliographystyle{amsplain}

\end{document}